\documentclass{tran-l}

\usepackage[active]{srcltx} 

\newtheorem{thm}{Theorem}
\newtheorem{lem}[thm]{Lemma}
\newcommand{\semi}{>\!\!\!\!\! \hspace{0.015cm}\lhd}
\newcommand{\To}{\longrightarrow}
\newcommand{\D}{\mathcal{D}}
\newcommand{\B}{\mathcal{B}}
\begin{document}

\title[Flag-transitive Steiner quadruple systems]
{Classification of flag-transitive \linebreak Steiner quadruple
systems}

\author{Michael Huber}

\address{Math. Institut der Universit\"{a}t T\"{u}bingen, Auf der Morgenstelle 10,
72076 T\"{u}bingen}

\email{mihu@michelangelo.mathematik.uni-tuebingen.de}






\begin{abstract}
A \emph{Steiner quadruple system} of order $v$ is a $3-(v,4,1)$
design, \linebreak and will be denoted $SQS(v)$.

Using the classification of finite $2$-transitive permutation
groups all $SQS(v)$ with a flag-transitive automorphism group are
completely classified, thus solving the "still open and
longstanding problem of classifying all flag-transitive
\mbox{$3-(v,k,1)$} designs"
(cf.~\cite[p.\,273]{Del95},\,\,\cite{Del92}) for the smallest
value of $k$. Moreover, a generalization of a result of H.
L\"{u}neburg \cite{Luene65} is achieved.
\end{abstract}

\maketitle

\section{Introduction}

In the last decades, there has been a great interest in
classifying $t-(v,k,\lambda)$ designs with certain transitivity
properties. For example, all point $2$-transitive \linebreak
\mbox{$2-(v,k,1)$ designs} were classified by Kantor~\cite{Kant85}
and a few years later Buekenhout et al.~\cite{Buek95} reached a
classification of all flag-transitive \mbox{$2-(v,k,1)$ designs.}
Both results depend on the classification of finite simple groups.
However, the classification of flag-transitive \mbox{$3-(v,k,1)$
designs} is "a still open and longstanding
problem"(cf.~\cite[p.\,273]{Del95},\,\,\cite{Del92}).

In this article we use the classification of finite $2$-transitive
permutation groups to classify all flag-transitive $SQS(v)$, thus
solving the above problem for the smallest value of $k$. Moreover,
our result generalizes a theorem of L\"{u}neburg~\cite{Luene65}
that characterizes all flag-transitive $SQS(v)$ under the
additional strong assumption that every non-identity element of
the automorphism group fixes at most two points. Our procedure as
well as our proofs are independent of L\"{u}neburg.

\smallskip

For positive integers $t \leq k \leq v$ and $\lambda$ we define a
\emph{$t-(v,k,\lambda)$ design} to be an incidence structure
$\mathcal{D}=(X,\mathcal{B},I)$, where $X$ is a set of
\emph{points}, $\left| X \right| =v$, and $\mathcal{B}$ a set of
\emph{blocks}, $\left| \mathcal{B} \right| =b$, with the
properties that each block $B \in \mathcal{B}$ is incident with
$k$ points, and every $t$-subset of $X$ is incident with $\lambda$
blocks. A \emph{Steiner quadruple system} of order $v$, which will
be denoted by $SQS(v)$, is a \mbox{$3-(v,4,1)$ design.}
Hanani~\cite{Han60} showed that a $SQS(v)$ exists if and only if
$v \equiv 2$ or $4$ (mod $6$) $(v \geq 4)$.

In the following let $\mathcal{D}=(X,\mathcal{B},I)$ be a
non-trivial $SQS(v)$ and $G \leq Aut (\mathcal{D})$ a group of
automorphisms of $\mathcal{D}$. A \emph{flag} is an incident
point-block pair, that is $x \in X$ and $B \in \mathcal{B}$ such
that $xIB$, and we call $G \leq Aut (\mathcal{D})$ to be
\emph{flag-transitive} (resp. \emph{block-transitive}) if $G$ acts
transitively on the flags (resp. on the blocks) of $\mathcal{D}$.
\linebreak For short, $\mathcal{D}$ is called
\emph{flag-transitive (resp. block-transitive, point
$t$-transitive)} if $\mathcal{D}$ admits a flag-transitive (resp.
block-transitive, point $t$-transitive) group of automorphisms.


Our result is the following

\begin{thm} \label{mainthm}
Let $\mathcal{D}=(X,\mathcal{B},I)$ be a non-trivial $SQS(v)$.
Then $G \leq Aut( \mathcal{D})$ acts flag-transitively on
$\mathcal{D}$ if and only if one of the following occurs:

\begin{enumerate}

\item[\emph{(1)}] $\mathcal{D}$ is isomorphic to the $SQS(2^d)$ whose
points and blocks are the points and planes of the affine space
$AG(d,2)$, and one of the following holds:
\item[] (i) $d \geq 3$, and $G \cong AGL(d,2)$,
\item[] (ii)  $d=3$, and $G \cong AGL(1,8)$ or $A \Gamma L(1,8)$,
\item[] (iii) $d=4$, and $G_0 \cong A_7$,
\item[] (iv) $d=5$, and $G \cong A \Gamma L(1,32)$,
\smallskip
\item[\emph{(2)}] $\mathcal{D}$ is isomorphic to a $SQS(3^d +1)$ whose
points are the elements of \linebreak \mbox{$GF(3^d) \cup
\{\infty\}$} and whose blocks are the images of $GF(3) \cup
\{\infty\}$ under \linebreak $PGL(2,3^d)$ with $d \geq 2$ (resp.
$PSL(2,3^d)$ with $d >1$ odd) and the derived design is isomorphic
to the \mbox{$2-(3^d,3,1)$ design} whose points and blocks are the
points and lines of $AG(d,3)$, and $PSL(2,3^d) \leq G \leq P
\Gamma L (2, 3^d)$,
\smallskip
\item[\emph{(3)}] $\mathcal{D}$ is isomorphic to a $SQS(q+1)$ whose
points are the elements of \linebreak \mbox{$GF(q) \cup
\{\infty\}$} with a prime power $q \equiv 7$ $(mod$ $12)$ and
whose blocks are the images of $\{0,1,\infty,\varepsilon\}$ under
$PSL(2,q)$, where $\varepsilon$ is a primitive sixth root of unity
in $GF(q)$ and the derived design is isomorphic to the Netto
triple system, \linebreak and $PSL(2,q) \leq G \leq P \Sigma L
(2,q)$.
\end{enumerate}
\end{thm}

\smallskip

A detailed description of the \emph{Netto triple system} can be
found in~\cite[Section\,3]{Deletal86}.
\bigskip \smallskip
\bigskip
\section{Preliminaries}

If $\mathcal{D}=(X,\mathcal{B},I)$ is a $t-(v,k,\lambda)$ design,
and $x \in X$ arbitrarily, the \emph{derived} design with respect
to $x$ is \mbox{$\mathcal{D}_x=(X_x,\mathcal{B}_x, I_x)$}, where
$X_x = X \backslash \{ x \}$, \mbox{$\mathcal{B}_x=\{B \in
\mathcal{B}: xIB \}$} and $I_x= I \!\!\mid _{X_x \times \;
\mathcal{B}_x}$. We shall also speak of $ \mathcal{D}$ as being an
\emph{extension} of $ \mathcal{D}_x$. \linebreak Obviously, a
derived design is a $(t-1)-(v-1,k-1,\lambda)$ design.

For $g \in G \leq Sym(X)$ let fix($g$) denote the set of fixed
points and supp($g$) the support of $g$. If $\{x_1,...,x_n\}
\subseteq X$ let $G_{\{x_1,...,x_n\}}$ be its setwise stabilizer
and $G_{x_1,...,x_n}$ its pointwise stabilizer. If $B \in
\mathcal{B}$ let $G_B$ be its block stabilizer and $G_{(B)}$ its
pointwise block stabilizer. By $r \perp q^n-1$ we mean that $r$
divides $q^n-1$ but not $q^k-1$ for all $1 \leq k < n$.

All other notation is standard.

\smallskip

If $\mathcal{D}=(X,\mathcal{B},I)$ is a $t-(v,k,1)$ design then it
is elementary that the point $2$-transitivity of $G \leq
Aut(\mathcal{D})$ implies its flag-transitivity when $t=2$.
However, for $t \geq 3$ the converse holds:

\begin{lem}\label{flag2trs}
Let $\mathcal{D}=(X,\mathcal{B},I)$ be a $t-(v,k,1)$ design with
$t \geq 3$. If $G \leq Aut( \mathcal{D})$ acts flag-transitively
on $ \mathcal{D}$ then $G$ also acts $2$-transitively on the
points of $\mathcal{D}$.
\end{lem}

\begin{proof}
Let $x \in X$. As $G$ acts flag-transitively on $ \mathcal{D}$,
obviously $G_x$ acts block-transitively on $\mathcal{D}_x$. Since
block-transitivity implies transitivity on points for $t \geq 2$
by Block's Theorem~\cite{Block65}, $G_x$ also acts transitively on
the points of $ \mathcal{D}_x$ and the claim follows.
\end{proof}

To classify all flag-transitive $SQS(v)$, we can therefore use the
classification of finite $2$-transitive permutation groups which
itself relies on the classification of finite simple groups
(cf.~\cite{CSK},\,\cite{Her74},\,\cite{Hup57},\,\cite{Kant85}).

The list of groups is as follows:

Let $G$ be a finite $2$-transitive permutation group of a
non-empty set $X$. \\Then we have either

{\bf (A) Affine type:} $G$ contains a regular normal subgroup $T$
which is elementary abelian of order $v=p^d$, where $p$ is a
prime. Let $a$ be a divisor of $d$. Identify $G$ with a group of
affine transformations
\[x \mapsto x^g+c\]
of $V(d,p)$, where $g \in G_0$. Then one the following occurs:
\begin{enumerate}
\item[(1)] $G \leq A \Gamma L(1,p^d)$
\item[(2)] $G_0 \unrhd SL(\frac{d}{a},p^a)$
\item[(3)] $G_0 \unrhd Sp(\frac{2d}{a},p^a)$, $d \geq 2a$
\item[(4)] $G_0 \unrhd G_2(2^a)'$, $d=6a$
\item[(5)] $G_0 \cong A_6$ or $A_7$, $v=2^4$
\item[(6)] $G_0 \unrhd SL(2,3)$ or $SL(2,5)$, $v=p^2$, $p=5,7,11,19,23,29$ or $59$, or $v=3^4$
\item[(7)] $G_0$ contains a normal extraspecial subgroup $E$ of order $2^5$, and $G_0/E$ is isomorphic
                 to a subgroup of $S_5$, where $v=3^4$
\item[(8)] $G_0 \cong SL(2,13)$, $v=3^6,$
\end{enumerate}

or \smallskip

{\bf (B) Semisimple type:} $G$ contains a simple normal subgroup
$N$, \linebreak and \mbox{$N \leq G \leq Aut(N)$}. In particular,
one of the following holds, where $N$ and $v=|X|$ are given:

\begin{enumerate}
\item[(1)] $A_v$, $v \geq 5$
\item[(2)] $PSL(d,q)$, $d \geq 2$, $v=\frac{q^d-1}{q-1}$, where $(d,q) \not= (2,2),(2,3)$
\item[(3)] $PSU(3,q^2)$, $v=q^3+1$, $q >2$
\item[(4)] $Sz(q)$, $v=q^2+1$, $q=2^{2e+1}>2$ \hfill (Suzuki group)
\item[(5)] $^2G_2(q)$, $v=q^3+1$, $q=3^{2e+1} > 3$ \hfill (Ree group)
\item[(6)] $Sp(2d,2)$, $d \geq 3$, $v = 2^{2d-1} \pm 2^{d-1}$
\item[(7)] $PSL(2,11)$, $v=11$
\item[(8)] $PSL(2,8)$, $v=28$ (N not $2$-transitive)
\item[(9)] $M_v$, $v=11,12,22,23,24$   \hfill (Mathieu group)
\item[(10)] $M_{11}$, $v=12$
\item[(11)] $A_7$, $v=15$
\item[(12)] $HS$, $v=176$ \hfill (Higman-Sims group)
\item[(13)] $Co_3$, $v=276$. \hfill (smallest Conway group)
\end{enumerate}

\bigskip

Let $r$ denote the number of blocks incident with a point. The
following obvious observation is important for this paper:
\begin{lem}\label{divprop}
Let $\mathcal{D}=(X,\mathcal{B},I)$ be a $t-(v,k,1)$ design, and
$x \in X$ arbitrarily. \linebreak If $G \leq Aut(\mathcal{D})$
acts flag-transitively on $\mathcal{D}$ then the division property
\[r \; \Big| \; \left| G_x \right|\] holds.
\end{lem}
Counting in two ways easily yields that $r=(v-1)(v-2)/6$ when
$\mathcal{D}$ is a $SQS(v)$.

\section{Proof of the theorem}

Using the notation as before, let $\D=(X,\B,I)$ be a $SQS(v)$. In
this section we run through the list of finite $2$-transitive
permutation groups given in Section $2$ and examine successively
whether $G \leq Aut(\D)$ acts flag-transitively on $\D$.

\subsection{Affine case}
From Section $2$ we know that a $2$-transitive permutation group
$G$ of affine type has degree $v=p^d$. As a $SQS(v)$ exists if and
only if \linebreak \mbox{$v=2$ or $4$ (mod $6$)} $(v \geq 4)$ by
Hanani's theorem, we conclude that $v=2^d$ in this case. To avoid
trivial $SQS(v)$, let $d \geq 3$.

The following lemma is fundamental for this case.
\begin{lem} \label{affSQS}
Let $\D=(X,\B,I)$ be a $SQS(2^d)$ with $d \geq 3$, and $G \leq
Aut(\D)$ contains a regular normal subgroup $T$ which is
elementary abelian of order $v=2^d$. \linebreak If $G$ acts
flag-transitively on $\D$ and $\left| G_0 \right| \equiv 1$ $(mod$
$2)$, then $\mathcal{D}$ is uniquely determined (up to
isomorphism), and the points and blocks of $\D$ are the points and
\linebreak planes of $AG(d,2)$.
\end{lem}

\begin{proof}
$T$ contains subgroups of order $4$ as it is elementary abelian of
order $2^d$. Moreover, $T$ is the only Sylow $2$-group since
$\left| G_0 \right| \equiv 1$ (mod $2$), and contains therefore
all subgroups of $G$ of order $4$. By assumption, $G_B$ acts
transitively on the points of $B$ for $B \in \B$ arbitrarily. Thus
$4$ is a divisor of the order of $G_B$, and $G_B$ contains at
least one subgroup $S$ of $T$ of order $4$. Then $B \in \B$ is an
orbit of $S$ and hence an affine plane. As $G \leq Aut(\D)$ is
block-transitive, we can conclude that all blocks must be affine
planes. Now identify the points of $\D$ with the elements of $T$
and the assertion follows.
\end{proof}

\emph{Case} (1): $G \leq A \Gamma L(1,2^d)$.

Let $\D=(X,\B,I)$ be a $SQS(2^d)$, $d \geq 3$, and assume $G \leq
Aut(\D)$ acts flag-transitively on $\D$. Lemma \ref{divprop} and
Lagrange's theorem yield
\[r=\textstyle{\frac{1}{3}}(2^d-1)(2^{d-1}-1) \Big| \left| G_0 \right| \Big|
\left| A \Gamma L(1,2^d)_0 \right| = \left| \Gamma L(1,2^d)
\right|=d(2^d-1).\] Thus $d=3,5$. First, assume $d=3$. Then
$\left| A \Gamma L(1,8) \right| = \left| T \right| \left| \Gamma
L(1,8) \right| = 8 \cdot 7 \cdot 3.$ Since $G$ is $2$-transitive,
we have $8 \cdot 7 \; \big| \left| G \right|$, hence $\left| G
\right|=8 \cdot 7$ or $8 \cdot 7 \cdot 3$. The latter implies $G
\cong A \Gamma L(1,8)$, so assume $\left| G \right|= 8 \cdot 7$.
Since $A \Gamma L(1,8)$ is solvable, we deduce from Hall's theorem
that $G \cong AGL(1,8)$ as $G$ is a Hall $\{2,7\}$-group. For
$d=5$ again $\left| G \right| =32 \cdot 31$ or $32 \cdot 31 \cdot
5$. We conclude $G \cong A \Gamma L(1,32)$ as for $\left| G
\right| = 32 \cdot 31$ lemma \ref{divprop} yields a contradiction.

On the contrary, we have to show that $G \cong AGL(1,8),A \Gamma
L(1,8)$ resp. $A \Gamma L(1,32)$ acts flag-transitively on the
$SQS(8)$ resp. the $SQS(32)$ given in the theorem. \linebreak For
$v=8$ there exists (up to isomorphism) only the unique $SQS(v)$
consisting of the points and planes of $AG(3,2)$. Since $G \cong
AGL(1,8)$ acts transitively on the points, it is sufficient to
show that $G_0 \cong GL(1,8)$ acts transitively on the blocks
incident with $0$. As these are exactly the $2$-dimensional
subspaces of the underlying vector space, we have \[B_1
:=\{0,1,t,t+1\} \neq B_1^{\,t}= \{0,t,t^2,t^2+1\} \quad \mbox{for}
\quad 1 \neq t \in GL(1,8) \cong GF(8)^*.\] Thus $\mid \!\!
B_1^{GL(1,8)}\!\! \mid \,\neq 1$, and hence as $r=7$, the claim
follows by the orbit-stabilizer property. Obviously, $G \cong A
\Gamma L(1,8)$ acts flag-transitively on $\D$ as well. For $v=32$
we have by lemma \ref{affSQS} also only the unique $SQS(v)$
consisting of the points and planes of $AG(5,2)$ because $\left|
G_0 \right| = \left| \Gamma L(1,32) \right| \equiv 1$ (mod $2$).
To see that \mbox{$G_0 \cong \Gamma L(1,32)$} acts
flag-transitively on the blocks incident with $0$, examine as
before that \mbox{$\mid \!\! B_1^{GL(1,32)}\!\! \mid \,\neq 1$},
thus \mbox{$\left| GL(1,32)_B \right|=1$} for any $0 \in B \in \B$
by the orbit-stabilizer property. Hence \mbox{$\mid \!\!B^{\Gamma
L(1,32)} \!\! \mid \, = 31$} or $31 \cdot 5$. Assuming the first
yields $\left| \Gamma L(1,32)_B \right| = 5$ by the
orbit-stabilizer property again. Let $H$ be a cyclic group of
order $5$. Then $\left| H_B \right| \neq 1$ for any $0 \in B \in
\B.$ On the other hand, $5$ is a $2$-primitive divisor of $2^4-1$.
Thus $H$ has irreducible modules of degree $4$ in view
of~\cite[Theorem\,3.5]{Her74}. As the $5$-dimensional
$GF(32)H$-module is completely reducible by Maschke's theorem, $H$
has as irreducible modules only the trivial module and one of
degree $4$. But if $H$ fixes any $2$-dimensional vector subspace
then, again by Maschke's theorem, $H$ would have as irreducible
modules two $1$-dimensional modules, a contradiction. Therefore,
$\mid \!\!B^{\Gamma L(1,32)} \!\!\mid \, = 31 \cdot 5$ must hold
and the claim follows as $r= 31 \cdot 5$.
\smallskip

\emph{Case} (2): $G_0 \unrhd SL(\frac{d}{a},2^a)$.

For $a=1$ we have $ G \cong AGL(d,2)$. Here $G$ is $3$-transitive
and the only $SQS(v)$ on which $G$ acts is the one whose points
and blocks are the points and planes of $AG(d,2)$, $d \geq 3$, by
Kantor~\cite{Kant85}. Obviously, $G$ is also flag-transitive. As
$a=d$ has already been done in \mbox{case (1)} we can assume that
$a$ is a proper divisor of $d$. We prove that here no
flag-transitive $SQS(v)$ exists.\\ Because of lemma~\ref{divprop},
it is enough to show that $r$ is no divisor of $\left| G_0
\right|$.\linebreak Clearly,
\[\left|SL(\textstyle{\frac{d}{a}},2^a)\right|=2^{d(\frac{d}{a}-1)/2}
\displaystyle{\prod_{i=2}^{\frac{d}{a}}(2^{ia}-1),}\] \[\mbox{and}
\quad [\Gamma L(\textstyle{\frac{d}{a}},2^a) :
SL\textstyle{(\frac{d}{a}} , 2^a)]=\left| Aut(GF(2^{a})) \right|
\left| GF(2^{a})^* \right| =a \cdot(2^a-1).\] Thus it is
sufficient to show that $r$ does not divide $a \cdot (2^a-1) \cdot
\left| SL(\textstyle{\frac{d}{a}},2^a) \right|$.\linebreak By
Zsigmondy's theorem (cf.~\cite[p.\,283]{Zsig})
\[2^{d-1}-1\]
has a  $2$-primitive prime divisor $\tilde{r} \neq 1$ with
$\tilde{r} \perp 2^{d-1}-1$. Obviously, $\tilde{r} \neq 2$.
Furthermore, $\tilde{r} \! \not \hspace{0.025cm} \!\!\;\,\! \mid
3a$ since $\tilde{r} \equiv 1$ $(mod$ $(d-1))$
(cf.~\cite[Theorem\,3.5] {Her74}) and $d$ is properly divisible by
$a$. \\ Therefore,
\[2^{d-1}-1 \not\!\!\;\big| \hspace{0.2cm} 3a \cdot 2^{d(\frac{d}{a}-1)/2}
\prod_{i=1}^{\frac{d}{a}-1} (2^{ia}-1)\] and the claim follows.
\smallskip

\emph{Cases} (3)-(4): These cases can be eliminated analogous case
(2) using lemma~\ref{divprop} and Zsigmondy's theorem. (For
$\left| Out(G_0) \right|$ see e.g.~\cite[Table \,5.1\,A]{KlLi}).
\smallskip

\emph{Case} (5): $G_0 \cong A_6$ or $A_7$, $v=2^4$.

If $G \cong A_6$ then lemma~\ref{divprop} implies that $G$ cannot
act flag-transitively on any $SQS(v)$.

As $ G \cong A_7$ is $3$-transitive and the only $SQS(v)$ on which
$G$ acts is the one whose points and blocks are the points and
planes of $AG(4,2)$ by Kantor~\cite{Kant85}, we have also
flag-transitivity in this case.
\smallskip

\emph{Cases} (6)-(8): These cases cannot occur since $v$ is no
power of $2$.


\subsection{Semisimple case}
The cases (3),\,(5),\,(8),\,(12) from the list where $G$ is of
semisimple type can easily be ruled out as above by using
lemma~\ref{divprop}. Obviously, \linebreak the cases
(4),\,(7),\,(10),\,(11),\,(13) cannot occur by Hanani's theorem.
\\ Before we proceed we indicate

\begin{lem} \label{PSL}
Let $V(d,q)$ be a vector space of dimension $d > 3$ over $GF(q)$
and $PG(d-1,q)$ the $(d-1)$-dimensional projective space. Assume
$G$ containing $PSL(d,q)$ acts on $PG(d-1,q)$ and for all $g \in
G$ with $\mid \!\! M^g \cap M \!\! \mid \;\geq 3$ we have $M^g =
M$, where $M$ is an arbitrary set of points of $PG(d-1,q)$ of
cardinality $k$ with $3 \leq k \leq \left| H \right|$, and $H$ a
hyperplane of $PG(d-1,q)$.\\ If $\left| M \cap H \right| \geq 3$,
then $M \cap H = M$ holds.
\end{lem}

\begin{proof}
For $k=3$ the assertion is trivial. So assume $3 < k \leq \left| H
\right| = \frac{q^{d-1}-1}{q-1}$. In $PG(d-1,q)$ Desargues'
theorem holds and the translations $T(H)$ form an abelian group
which is sharply transitive on the points of $PG(d-1,q) \setminus
H$ by Baer's theorem. But on $H$ the group $T(H)$ acts trivially
since the central collineations fix each point of $H$. Thus the
claim holds if all elements of $M$ lie in $H$. Therefore, assume
that there is an element of $M$ which is not in $H$. Then $M$
contains all points of $PG(d-1,q) \setminus H$ as $T(H)$ is
transitive. Thus
\[\left| M \right| \geq \frac{q^d-1}{q-1} - \frac{q^{d-1}-1}{q-1}=
\frac{q^d-q^{d-1}}{q-1}=q^{d-1} > \frac{q^{d-1}-1}{q-1}=\left| H
\right|.\] But this contradicts the assumption $\left| M \right|
\leq \left| H \right|,$ and the claim follows.
\end{proof}

\emph{Case} (1): $N=A_v$, $v \geq 5$. Here, $G$ is $3$-transitive
and does not act on any non-trivial $3-(v,k,1)$ design by
Kantor~\cite{Kant85}. \smallskip

\emph{Case} (2): $N=PSL(d,q)$, $d \geq 2$, $v=\frac{q^d-1}{q-1}$,
where $(d,q) \not= (2,2),(2,3)$.

We distinguish two subcases:

(i) $N=PSL(2,q)$, $v=q+1$.

Here $q \geq 5$ as $PSL(2,4) \cong PSL(2,5)$, and $Aut(N)= P
\Gamma L (2,q)$. First suppose that $G$ is $3$-transitive.
According to Kantor~\cite{Kant85}, we have then only the
$SQS(3^d+1)$ described in (2) of theorem~\ref{mainthm}, and
$PSL(2,3^d) \leq G \leq P \Gamma L(2,3^d)$. Obviously, also
flag-transitivity holds. As $PGL(2,q)$ is a transitive extension
of $AGL(1,q)$, it is easily seen that the derived design at any
point of $GF(3^d) \cup \{\infty\}$ is isomorphic to the
$2-(3^d,3,1)$ design consisting of the points and lines of
$AG(d,3)$.

Now assume that $G$ is $3$-homogeneous but not $3$-transitive. As
here $PSL(2,q)$ is a transitive extension of $AG^2L(1,q)$ we
deduce from~\cite{Deletal86} that the derived design is either the
affine space $AG(d,3)$ or the Netto triple system. Thus (2) with
the part in brackets or (3) of theorem~\ref{mainthm} holds with
$PSL(2,3^d) \leq G \leq P\Sigma L(2,3^d)$ (where \linebreak $P
\Sigma L (2,p^d):= PSL(2,p^d) \semi <\!\!\tau_{\alpha} \!\!>$ with
$\tau_{\alpha} \in Sym(GF(p^d) \cup \{\infty\}) \cong S_v$ of
order $d$ induced by the Frobenius automorphism $\alpha : GF(p^d)
\longrightarrow GF(p^d),\, x \mapsto x^p$.) Conversely, as $G$ is
$3$-homogeneous it is also block-transitive. In both cases we have
$PSL(2,q)_B \cong A_4$ for any $B \in \B$ since $PSL(2,q)_B$ has
order $12$ by the orbit-stabilizer property and $PSL(2,q)_B
\longrightarrow Sym(B) \cong S_4$ is a faithful representation.
Thus, in each case flag-transitivity holds.

Finally, suppose $G$ is not $3$-homogeneous. As $PGL(2,q)$ is
$3$-homogeneous the $PGL(2,q)$-orbit on $3$-subsets therefore
splits under $PSL(2,q)$ into two orbits of same length. Let $M$ be
an arbitrary $3$-subset. Then $\left| PSL(2,q)_M \right| = \left|
PGL(2,q)_M \right| = 6$ by the orbit-stabilizer property.
\pagebreak Thus, as $PGL(2,q)$ is \mbox{$3$-transitive} we have
$PSL(2,q)_M \cong S_3$ for each orbit. If $PSL(2,q)$ acts
block-transitively on any $SQS(v)$ then $PSL(2,q)_B \cong A_4$
again for any $B \in \B.$ But, by the definition of $SQS(v)$ this
would imply that $PSL(2,q)_{\tilde B}$, where $\tilde{B}$ denotes
the block uniquely determined by $M$, contains $PSL(2,q)_M$, a
contradiction. Thus $PSL(2,q)$ does not act flag-transitively on
any $SQS(v)$. We show now that $G$ cannot act flag-transitively on
any $SQS(v)$. Without restriction choose $\mathcal{O}_1$ to be the
\mbox{$PSL(2,q)$-orbit} containing $\{0,1,\infty\}$. Easy
calculation shows that \mbox{$P \Sigma L (2,q)_{0,1,\infty} = \;<
\!\! \tau_{\alpha} \!\!>$}. Thus $P \Sigma L(2,q)_{\mathcal{O}_1}$
is contained in $P \Gamma L(2,q)$, and equality holds as $P \Sigma
L(2,q)$ is of \mbox{index $2$} in $P\Gamma L(2,q)$ and $P \Gamma
L(2,q)$ is $3$-transitive. Therefore, we only have to consider
\mbox{$PSL(2,q) \leq G \leq P \Sigma L(2,q)$.} Dedekind's law
yields \mbox{$G = PSL(2,q)  > \!\!\! \lhd \;(G \;\,\cap < \!\!
\tau_{\alpha} \!\!>)$} and \mbox{$G_{(B)} = PSL(2,q)_{(B)}  >
\!\!\! \lhd \; G \;\,\cap < \!\! \tau_{\alpha} \!\!> \; = \,G
\;\,\cap < \!\! \tau_{\alpha} \!\!> \,\cong C_m$}, the cyclic
group of order $m \,| \, d$, for any $B \in \B$ since every
non-identity element of $PSL(2,q)$ fixes at most two points.
Assume $G$ acts block-transitively on any $SQS(v)$. Then we can
choose $B \in \B$ such that $B$ contains $\{0,1,\infty\}$. Since
$G_{(B)}$ is the kernel of the representation $G_B \To Sym(B)
\cong S_4$ and $PSL(2,q)_B \cong A_4$ we have therefore again by
Dedekind's law
\[G_B = PSL(2,q)_B \times (G \;\,\cap < \!\! \tau_{\alpha} \!\!>) \cong A_4 \times C_m.\]
However, as $PSL(2,q)_{\{0,1,\infty \}} \cong S_3$ we get
analogously
\[G_{\{0,1,\infty \}} = PSL(2,q)_{\{0,1,\infty \}} \times (G \;\,\cap < \!\! \tau_{\alpha} \!\!>) \cong S_3 \times
C_m,\] which leads again to a contradiction by the definition of
$SQS(v)$.

(ii) $N=PSL(d,q)$, $d \geq 3$, $v=\frac{q^d-1}{q-1}.$

Here $Aut(N)=P \Gamma L(d,q) > \!\!\!\! \lhd  <\!\! \iota \!\!>$,
where $\iota$ denotes a graph automorphism. We show that $G$ does
not act on any $SQS(v)$. For $d=3$ this is obvious since
$v=q^2+q+1$ is always odd, a contradiction to Hanani's theorem.

Consider $d>3$ and let $H$ be a hyperplane of the projective space
$PG(d-1,q)$. Assume that the claim does not hold. Then there is a
counterexample with $d$ minimal. Without restriction we can choose
three arbitrary points $\alpha ,\beta, \gamma$ from $H$. As for
$d>3$
\[|H|=\frac{q^{d-1}-1}{q-1} >4\]
holds, the block uniquely determined by $\alpha ,\beta, \gamma$ is
contained in $H$ by lemma \ref{PSL}. Thus $H$ induces a
$SQS(\frac{q^{d-1}-1}{q-1})$ on which $G$ containing $PSL(d-1,q)$
operates. By induction, we get the minimal counterexample for
$d=3$. So $G$ containing $PSL(3,q)$ acts on a
$SQS(\frac{q^3-1}{q-1})$. But, as above
$\frac{q^3-1}{q-1}=q^2+q+1$ is always odd yielding the desired
contradiction.

\smallskip

\emph{Case} (6): $N=Sp(2d,2)$, $d \geq 3$, $v = 2^{2d-1} \pm
2^{d-1}$.

Here $N=G$ since \mbox{$\left|Out(N) \right|=1$} (cf.~\cite[Table
\,5.1\,A]{KlLi}). We show that $G$ contains elements which fix
exactly $3$ points and hence cannot act on any $SQS(v)$ by
definition.

Let $X^+$ respectively $X^-$ denote the set of points on which $G$
operates with \linebreak \mbox{$\left| X^+ \right|=2^{2d-1} +
2^{d-1}$} resp. \mbox{$\left| X^- \right|=2^{2d-1} - 2^{d-1}$} ,
and define
\[m_p(G) := \mbox{min} \{ \left| \mbox{supp}(g) \right|  : 1 \neq
g \in G, g \,\mbox{ a }\, p\mbox{-element of }G \}\] to be the
\emph{minimal $p$-degree} of a transitive permutation group $G$,
$p$ a prime divisor of $\left| G \right|$ (cf.~\cite{Jens}).

First, suppose $d$ is even. By Zsigmondy's theorem
\[2^{d-1}-1\]
has a $2$-primitive prime divisor $p$ with $p \perp 2^{d-1}-1$.
Moreover, $p$ divides $\left| G \right|$ since $\left| G \right| =
2^{d^2} \prod^d_{i=1} (2^{2i}-1)$ (see e.g.~\cite[Table
\,2.1\,C]{KlLi}). Therefore, according to~\mbox{\cite[Theorem
3.7]{Jens}} we get in $X^+$
\[m_p(G) = 2^{2d-2(d-1)-1} (2^{2(d-1)}- 1) +
2^{d-(d-1)-1}(2^{d-1}-1)= \left| X^+ \right| - 3.\] Thus, there
exists $g \in G$ of prime order $p$ that fixes 3 points in $X^+$.

For $d \neq 4$ Zsigmondy's theorem yields the existence of a
$2$-primitive prime divisor $p$ with $p \perp 2^{2(d-1)}-1$ and as
$p$ divides $\left| G \right|$ we have in $X^-$ again
by~\mbox{\cite[Theorem 3.7]{Jens}}
\[m_p(G) = 2^{2d-2(d-1)-1} (2^{2(d-1)}- 1) - 2^{d-(d-1)-1}(2^{d-1}+1)= \left|
X^- \right| - 3.\] When $d=4$ then~\cite[p.\,123]{Atlas} yields
$\left| \mbox{fix}(g) \right|=3$ in $X^-$ for $g \in 3D$, where
$3D$ denotes a conjugacy class in~\cite{Atlas}.

Now, suppose $d$ is odd. Again by Zsigmondy's theorem
and~\cite[Theorem 3.7]{Jens} there exists a $2$-primitive prime
divisor $p$ with $p \perp 2^{2(d-1)}-1$, and $m_p(G)= \left| X^-
\right| -3$ in $X^-$.

If $d \neq 7$ Zsigmondy's theorem yields the existence of a
$2$-primitive prime divisor $p$ with $p \perp 2^{d-1}-1$. Choose
$\begin{pmatrix}
  A_0 & A_1 \\
  A_2 & A_3
\end{pmatrix} \in S \in Syl_p(Sp(d-1,2))$ and define

\[h:=\begin{pmatrix}
  A_0 &  &  & A_1 &  &  \\
   & A_0 &  &  & A_1 &  \\
     &  & 1 &  &  & 0 \\
  A_2 &  &  & A_3 &  &  \\
   & A_2 &  &  & A_3 &  \\
   &  & 0 &  &  & 1 \\
\end{pmatrix}.\] The proof of~\cite[Theorem 3.7]{Jens} yields \mbox{$\left| \mbox{fix}(h)
\right|=3$} in $X^+$ and \mbox{$\left| \mbox{fix}(h) \right|=1$}
in $X^-$. \\ For $d = 7$ choose $A:=\begin{pmatrix}
  1 & 1 \\
  1 & 0
\end{pmatrix}$ and define \mbox{${k}:=diag(A,\,A,\,A,\,1,\,^t
\!\!A^{-1},\,^t\!\!A^{-1},\,^t\!\!A^{-1},\,1)$.} Again,
\mbox{$\left| \mbox{fix}({k}) \right|=3$} in $X^+$ and
\mbox{$\left| \mbox{fix}({k}) \right|=1$} in $X^-$. Thus the
assertion is proved.

\smallskip \smallskip

\emph{Case} (9): $M_v$, $v=11,12,22,23,24$.

Here, only $v=22$ is possible by Hanani's theorem. But as $M_{22}$
is $3$-transitive, Kantor~\cite{Kant85} shows that the only
$3-(v,k,1)$ design on which $M_{22}$ resp. $Aut(M_{22})$ acts is
the $3-(22,6,1)$ design. Therefore, this case cannot occur
finishing the proof of theorem~\ref{mainthm}.

\subsection*{Acknowledgment}
I would like to thank my supervisor Ch. Hering for his advice and
helpful discussions.

\INPUT{Xbib.bib}   
\INPUT{Paper.bbl}  
\bibliographystyle{amsplain}
\bibliography{xbib}
\end{document}